\documentclass[11pt]{article}

\usepackage{fullpage}
\usepackage{epsfig}
\usepackage{algorithm}
\usepackage{algorithmic}
\usepackage{amsmath}
\usepackage{amssymb}
\usepackage{amsfonts}
\usepackage{psfrag}
\usepackage{bbm}
\usepackage{listings}
\usepackage{cite}
\usepackage{array}

\newtheorem{definition}{Definition}

\title{A Dynamic Algorithm for Facilitated Charging of \\ Plug-In Electric Vehicles}
\author{Nicole Taheri%
   \thanks{Institute for Computational and Mathematical Engineering,
           Stanford University, Stanford, CA.
           Email: {\tt ntaheri@stanford.edu}.} 
   \and Robert Entriken\thanks{Electric Power Research Institute, Palo 
	 	Alto, CA. Email: {\tt rentrike@epri.com}.}
   \and Yinyu Ye%
   \thanks{Department of Management Science and Engineering,
           Stanford University, Stanford, CA.
           Email: {\tt yinyu-ye@stanford.edu}.}}
\date{\today}

\newcommand{\url}[1]{ {\tt #1 }}

\newcommand{\R}{\mathbbm{R}}

\newcommand{\matlab}{{\sc Matlab}}

\newcommand{\bmat}[1]{\begin{bmatrix}#1\end{bmatrix}}

\newcommand{\compresslist}{
\setlength{\itemsep}{0pt}
\setlength{\parskip}{0pt}
\setlength{\parsep}{0pt}
\setlength{\rightmargin}{10em}}

\newcommand{\minimize}{\text{minimize}}
\newcommand{\maximize}{\text{maximize}}
\newcommand{\st}{\text{subject to}}

\newcommand{\kpm}{\text{kWh/mile}}
\newcommand{\eff}{\text{eff}}

\begin{document} 
\maketitle
\begin{abstract}
Plug-in Electric Vehicles (PEVs) are a rapidly developing technology that can
reduce greenhouse gas emissions and change the way vehicles obtain power.
PEV charging stations will most likely be available at home and at work, and
occasionally be publicly available, offering flexible charging options.
Ideally, each vehicle will charge during periods when electricity prices are
relatively low, to minimize the cost to the consumer and maximize societal
benefits.  A Demand Response (DR) service for a fleet of PEVs could yield
such charging schedules by regulating consumer electricity use during certain
time periods, in order to meet an obligation to the market.  

We construct an automated DR mechanism for a fleet of PEVs that facilitates
vehicle charging to ensure the demands of the vehicles and the market are
met.  Our dynamic algorithm depends only on the knowledge of a few hundred
driving behaviors from a previous similar day, and uses a simple adjusted
pricing scheme to instantly assign feasible and satisfactory charging
schedules to thousands of vehicles in a fleet as they plug-in. The charging
schedules generated using our adjusted pricing scheme can ensure that a new
demand peak is not created and can reduce the consumer cost by over 30\% when
compared to standard charging, which may also increase peak demand by 3.5\%.
In this paper, we present our formulation, algorithm and results.
\end{abstract}

\section{Introduction}
A Plug-in Electric Vehicle (PEV) is any vehicle that uses electricity from the
grid to displace liquid fuel. These vehicles will be able to plug-in to the
grid and control which times during the connection period the vehicle battery
will actually charge.  This capability will allow consumers to charge their
vehicles during off-peak hours,  usually also resulting in a reduced cost to
the consumer.  For example, suppose a vehicle is plugged in from 6pm until 7am,
giving an 11-hour window in which the vehicle can charge.  If the battery can
obtain sufficient charge in one hour, this charge can be scheduled to occur
during an ideal hour between 6pm and 7am, which takes into account the
electricity price, electricity load, and individual vehicle battery state of
charge and driving schedule.  This will reduce consumer cost of charging and help
cultivate a balance between the electricity supply and the demand.

In order for PEV charging to be advantageous for both the PEV owners and the
utility, a mechanism is needed that receives information from both the vehicles
and the system dispatcher to facilitate vehicle charging. In this paper, we
design an algorithm for such a mechanism that will be used by an aggregator.

\begin{definition}
	An \emph{aggregator} manages the communication and electricity distribution
	between a group of electricity consumers and an electric utility. 
\end{definition}

Currently, aggregators regulate consumer electricity use by rewarding decreases
in electricity demand during certain time periods in order to meet a scheduling
obligation to the market or system dispatcher.  These Demand Response (DR)
resources are provided by either third-parties (such as EnerNOC and ZigBee
Smart Energy) or the utilities themselves (for example, both PG\&E and Duke
Energy have DR programs, among many others).  Existing programs focus on the
Commercial and Industrial Sector, with a small number of large customers who
typically use hundreds of kilowatt hours (kWh) per day \cite{EIA-usage}.  These
programs are not automated, requiring human intervention to individually tailor
the DR resource for each customer and monitor the electricity use of each
customer alongside the market supply and prices. 

By contrast, some estimates say there may be 100 million PEVs on the road in
the United States by 2030 \cite{epri-prism} and each PEV battery will generally
require a few kWh each day.  Thus, an automated DR program would be more
convenient and scalable to large fleets of vehicles than the alternative of
individually tailoring human-operated DR for each consumer.  An automated DR
service can manage the transactions between a fleet of PEVs (say, around
10,000) and their utility, facilitating the charging schedules of the fleet to
meet the transportation demands of the vehicle owners and meet a scheduling
obligation. 

In this paper, we describe an automated mechanism that performs Demand Response
services with a fleet of Plug-In Electric Vehicles. Our algorithm constructs
adjusted `prices' for each hour that take into account a scheduling obligation,
transportation demands of the fleet, vehicle driving schedules, and current
electricity pricing.  Vehicle charging schedules will be based on these
adjusted prices and determined instantly as vehicles plug-in, without knowing
the energy needs of subsequent vehicles connecting.  The resulting charging
schedules will meet the energy demands of all vehicles in the fleet while
avoiding increased load during certain time periods, in order to meet the
scheduling obligations.  For example, our algorithm can help ensure that a new
peak is not created while meeting vehicle charging needs, or that the added amount of
electrical load at any time is below a given threshold, which is a key feature
to maintaining the stability of the power grid.  

Our mechanism can be implemented as a device attached to each vehicle in a
fleet; these devices can communicate with each other and regulate vehicle
charging.  Installing the device on a vehicle will result in two main benefits
for the consumer: First, the cost of charging will be reduced since the vehicle
will charge in feasible hours when the electricity cost is lowest.  Second, the
device will decrease the risk of insufficient energy in the battery to drive;
as opposed to a driver simply plugging in a PEV arbitrarily, the driver will be
notified of the best times to plug-in to ensure the vehicle is sufficiently
charged.  The utilities will also benefit from a DR mechanism, since it will
maintain the same peak power demand, whereas PEV charging without demand
response will increase the peak power demand.  Such a device is a practical
extension of services that are currently provided with PEV charging stations,
which allow the user to control charging but require settings to be determined
by the user \cite{blink}.  

An aggregator's role is between these devices and the dispatcher to establish
and monitor the market supply and vehicle demand.  We assume such an
aggregator will exist, and that it will be possible for users to communicate
to the aggregator their expected driving schedule and vehicle characteristics.
It should be mentioned that such managed charging is often envisioned to be an
`opt in' service, in which a consumer has the option to participate.

\subsection{Motivation}
\label{motivation}

It is projected that by the year 2030, between 6\% and 30\% of vehicles in use
will be PEVs \cite{impact-EPRI}.  These electric vehicles will each run on
energy provided by a battery that will charge from the electricity grid, which
suggests an increase in electricity demand.  The Electric Power Research
Institute (EPRI) showed in their work \cite{impact-EPRI} that in the worst
case, the increase in total grid resource capacity will be 5-6\%, and \emph{smart
charging}, i.e., shifting some charging to off-peak hours, will decrease the
impact to only 1-2\%.  In a collaborative work between Better Place and PJM
Interconnection, Schneider et al.\ showed that such controlled charging will
reduce consumer energy costs by 45\% \cite{pjm-bp}. 

As more motivation, renewable energy sources are becoming increasingly
prevalent.  Many renewable power plants generate more electricity than is
needed during off-peak hours, and a lack of sufficient energy storage means
this surplus cannot always be stored for later use to meet the demand.  For
example, during the first half of 2008, an overproduction of wind power in
West Texas led to negative electricity prices 20\% of the time, and most of
this surplus occurred in the middle of the night \cite{TX-neg}. If vehicles
had been connected to the grid with a DR mechanism, they may have actually
profited from charging their batteries during these time periods.  

\subsection{Related Work}

A number of papers in the field of electric transportation have established the
benefits of smart charging, including \cite{impact-EPRI,irc10,nrel10,
phevs-vermont}.  However, there is currently no standard agreement on how to
manage PEV charging.  No previously proposed algorithm is based on a relatively
small amount of information and dynamically creates charging schedules for
thousands of vehicles in a fleet while meeting an external scheduling
obligation.

Han et al.\ in \cite{v2g-DP10} use dynamic programming to assign charging
schedules that provide frequency regulation. Similarly, in the work by Wu et
al.\ \cite{wu11} an algorithm is constructed that makes dynamic decisions for
lowest-cost charging schedules of PEVs.  Neither of these works considers an
obligation to the market that needs to be met, and the demand of each vehicle
is considered individually (i.e., not as part of a fleet); both suggested
algorithms would result in an increased peak demand.  In \cite{decent10}, Ma et
al.\ establish a decentralized algorithm that determines an equilibrium price
so that the total amount of charging done in the fleet fills the `overnight
demand valley.'  This algorithm takes into account all vehicles in the fleet,
but is not dynamic (i.e., it is solves an offline problem) and assumes only a
few types of driving behaviors exist.  A dynamic (or online) algorithm does not
require that all vehicles are connected to the grid at the same time to exactly
report their future driving schedules, and makes decisions for each vehicle
without solving a problem that depends on knowing the demands of the entire
fleet.

Our algorithm depends only on historical data, or a few hundred driving
behaviors from a previous similar day.  The implementation is very simple: a
set of unique adjusted `prices' are generated for each driving behavior and
charging is allocated among hours with the lowest adjusted price value.  This
pricing scheme will distribute electricity to thousands of PEVs in a fleet
while meeting an obligation to the market.  It is assumed that vehicles will
connect to the grid during periods when they are parked for at least an hour.
We use real data on driving behaviors, electricity loads, electricity and
gasoline pricing, and vehicle characteristics to generate our results. 

\section{Problem Description}

An aggregator will make commitments to both the PEV consumers and their
utility.  To the utility, the aggregator will commit to ensuring the
electricity supplied to the vehicles stays within a certain upper limit; we say
that the market places a \emph{charge cap} on the total amount of charging that
can be done by the fleet in each hour, in order to stay within the given upper
limit.  In addition, vehicles will plug-in to the grid at various times and the
aggregator will commit to providing each PEV with enough energy to drive.
Thus, each vehicle's charging decision should:
\begin{enumerate}\compresslist
	\item be determined instantly as vehicles plug-in to the grid,
	\item provide enough charge to each individual PEV to meet its daily
		transport load,
	\item meet the obligation to the market,
	\item consider the demand of the fleet as a whole, and
	\item reduce the cost to the consumer.
\end{enumerate}

Our algorithm is based on linear programming theory and achieves all of the
above for a time period of $n$ hours in the future; the value of $n$ is easily
adjustable.  For instance, our results consider $n=120$ hours (or 5 days) in
the future and solve the optimization problem \eqref{lp_clst} once at the start
of the 5-day period to determine charging decisions for the next 120 hours.
Vehicles are assumed to plug-in periodically over a given period of time (say,
12 hours) and report their driving schedule for the next $n$ hours.  Moreover,
we assume the information on the vehicle specifications can be obtained and
that the aggregator has collaborated with the market to place a charge cap on
the total amount of charging that can be done in each hour by the fleet.  In
the dynamic algorithm, each vehicle's charging schedule is instantly determined
as it connects to the grid. We assume a few hundred (say, 400) driver behaviors
from a similar day in a similar region are known, in addition to the predicted
electricity load for the next $n$ hours. 

\subsection{Linear Programming Formulation}

For explanation purposes, in the following linear programming formulation we
assume perfect foresight so that the aggregator has all the information ahead
of time (this is not assumed for our dynamic algorithm).  Let $m$ be the number
of vehicles and $n$ be the number of hours in the future considered.  We write
the driving load of each vehicle $j$ as a vector $t_j \in \R^{n}$, where the
$h$th element $t_{jh}$ is the number of miles driven in hour $h$.  The
aggregator also knows the electricity prices at each hour $p \in \R^n$, the
gasoline prices $p^g \in \R^n$, the market's charge cap for each hour
$c_{\text{cap}} \in \R^n$, and the characteristics of each vehicle.  Our
mechanism considers Plug-in Hybrid Electric Vehicles (PHEVs) that receive power
from both gasoline and electricity, in addition to Battery Electric Vehicles
(BEVs) that only receive power from electricity.  A full list of parameters and
variables is given in Table \ref{var-list}.

\begin{table}[!htbp]
\begin{center}
	\small
	\begin{tabular}{|c|l|c|l|c|}
		\hline
		\multicolumn{2}{|c|}{Variables} & \multicolumn{2}{|c|}{Input
			Parameters} & Dimension  \\
		\hline
		$s_j$ & electricity storage  & $\bar{s}_j$ & storage
			capacity  & $ n$ \\
		\hline
		$c_j$ & charging schedule & $\bar{c}_j$ & maximum charging
			rate  & $ n$ \\
		\hline
		$s^g_j$ & gasoline storage & $\bar{s^g}_j$ & gas tank
			capacity  & $ n$\\
		\hline
		$f_j$ & fueling schedule & $\bar{f}_j$ & maximum fueling
			rate & $ n$ \\
		\hline
		$g_j$ & generation schedule & $\bar{g}_j$ & maximum
			generation rate  & $ n$\\
		\hline
		&& $t_j$ & driver transport load &  $ n$\\
		\hline
		&& $p,p^g$ & electricity and gasoline prices & $ n$\\
		\hline
		&& $c_\text{cap}$ & charge cap	& $n$ \\
		\hline
		&& $s_{j0}, s^g_{j0}$ & initial battery and gasoline storage
			& $1$ \\
		\hline
		&& c$_\eff$, g$_\eff$ & charging and generating efficiencies
			& $1$ \\
		\hline
	\end{tabular}
	\caption{List of Parameters and Variables}
	\label{var-list}
\end{center}
\end{table}
The following linear program \eqref{lp} will find feasible charging, fueling,
and generating schedules, $c_j,f_j,g_j \in \R^{n}$, for each vehicle
$j=1,2,\ldots,m$, such that each vehicle has enough energy in its battery to
drive, while minimizing the cost and staying within the charge cap.  It is
assumed here that the aggregator has perfect foresight for the $n$-hour period.
We use \matlab \ vector notation, where $x_{a:b}$ refers to the elements
$\{a,a+1, \ldots, b\}$ of $x$.

\begin{equation}
	\begin{array}{lll}
		\underset{c_j, f_j, g_j, \forall j}{\minimize} &
			\sum_{j=1}^m \left(p^T c_j + (p^g)^Tf_j \right) \\
		\st & s_{j,1:n} = s_{j,0:n-1} + c_{\text{eff}} \cdot c_j + 
			g_{\text{eff}} \cdot g_j - (\kpm) \cdot t_j \\
	 	& s^g_{j,1:n} = s^g_{j,0:n-1} + f_j - (\text{gallon/kWh})
		 	\cdot g_j \\
		& t_{jh} = 0 \; \implies \; f_{jh} = 0, g_{jh} = 0 \\
		& t_{jh} > 0 \; \implies \; c_{jh} = 0 \\
		& \sum_{j=1}^m c_j \le c_{\text{cap}}\\
		& 0 \le (s_j,c_j,s^g_j,f_j,g_j) \le (\bar{s}_j,\bar{c}_j,
			\bar{s}^g_j,\bar{f}_j,\bar{g}_j) 
	\end{array}
	\label{lp}
	\tag{LP}
\end{equation}
The first constraint updates the battery storage amount and the second updates
the gasoline storage amount.  The third and fourth constraints ensure that the
vehicle can only fuel the tank and generate electricity when the vehicle is
driving, and that the vehicle battery can only be charged when the vehicle is
parked for at least an hour.  The second to last constraint ensures that the
total amount of charging done in each hour is below the charge cap,
$c_{\text{cap}}$.  And the last constraint ensures the variables stay within
their physical bounds.  

This problem formulation is static, since it assumes the aggregator knows the
exact future driving schedules of all vehicles that plug-in to the grid.
Moreover, our testing showed that solving this problem for a fleet of 10,000
vehicles over a 5-day period takes hours on a single workstation.  Thus,
instead of considering the micro-level linear program described above for all
vehicles, we only consider a few hundred possible vehicle behaviors and their
aggregated demands in this linear program.  These data can be based on previous
driving schedules from a relatively small population.  After solving the
aggregated linear program, we obtain adjusted `prices' and use them to
\emph{dynamically} assign feasible and satisfactory charging schedules to
thousands of individual vehicles within milliseconds of a connection to the
grid.

\section{Clustering}

Our algorithm assigns charging schedules that consider possible driving
behaviors of an entire fleet.  Since it is impossible to predict the exact
driving behaviors and available charging times of individual vehicles of a
fleet, we use clustering to estimate the expected future demand. This
approximation of future driving loads is necessary to stay below a given charge
cap.

\subsection{Lowest-Cost Charging}

Allocating charge to each vehicle in the hours with the lowest electricity price
will guarantee the lowest cost to each driver, but generally results in an
increase in peak electricity demand.  This is the method used in \cite{wu11},
and is not feasible when considering scheduling obligations.  For example, consider
the following scenario with two pure electric vehicles in a three hour period,
where
\[
p = \bmat{0.10 \\ 0.12 \\ 0.14}, \quad 
t_1 = \bmat{ 0 \\ 0 \\ 1 }, \quad
t_2 = \bmat{ 0 \\ 1 \\ 0 }, \quad
c_\text{cap} = (\kpm) \cdot \bmat{ 1 \\ 1 \\ 1},
\quad \text{and} \quad
s_{10} = s_{20} = 0.
\]

Assume both vehicles connect to the grid before hour 1, but at slightly
different times, and charging schedules are allocated on a first-come,
first-served basis.  If vehicle 1 arrives first, an assignment based on
electricity price will charge vehicle $1$ in hour 1, and the given charge cap
implies there will be no more electricity available in hour 1.  When vehicle 2
plugs in slightly later, it will not be able to charge during hour 1 since the
charge cap is already reached.  Thus, vehicle 2 will not have sufficient
battery energy storage to drive, and the algorithm has failed by not
considering market constraints and the demand of future vehicles to connect.
Note that because we assume a model where vehicles pay a low-cost flat rate for
facilitated charging, vehicle 1 will not end up paying more than vehicle 2 and
vice versa. 

Note that for simplicity in explanation, we say in this example that the charge
cap is in terms of kWh (energy), not kW (power), however, this is equivalent to
saying the charge cap is a constant 1 kW of power, which totals to 1 kWh of
energy over each hour. 

\subsection{Clustering Driving Patterns}

In order to consider future demands of vehicles in a fleet, a clustering
algorithm is used to determine a number of different \emph{base} driving
profiles.  We define a \emph{base} driving profile to be the ``best-fit'' to a
group of drivers with similar driving patterns.  Each vehicle is then matched
with the most similar base driving profile and properties of the corresponding
cluster are used to determine the individual charging schedule.  

Considering the approximate needs of each cluster before assigning individual
schedules will ensure the demands of each vehicle can be met.  Moreover, the
aggregator can use the base driving profiles and the expected number of drivers
in each cluster to approximate the energy demand of the fleet, which is useful
in negotiations with the market on placing a charge cap (i.e., in determining
scheduling obligations).

Real driving behavior data was used for the clustering and was obtained from
the National Household Transportation Survey (NHTS) \cite{nhts}. The survey
contains the driving schedule of over 150,000 individuals, each for a 24-hour
period.  We took data for each individual and put it into vector format, so
each individual vehicle $j=1,\ldots,m$ has a vector of the form:
\[
t_j = \bmat{ t_{j1}, t_{j2}, \cdots, t_{j,24} } \in \R^{24}.
\]
For simulation purposes, these transport loads were assigned to either a 
PHEV or a BEV: if the daily transport load amounted to less than 70 miles, the
corresponding vehicle was assumed to be a BEV; otherwise, it was assumed to be
a PHEV.  Realistically, the vehicle type would be information passed from the
drivers to the aggregator.
\begin{figure}[!htbp]
\begin{center}
	\includegraphics[scale=.65]{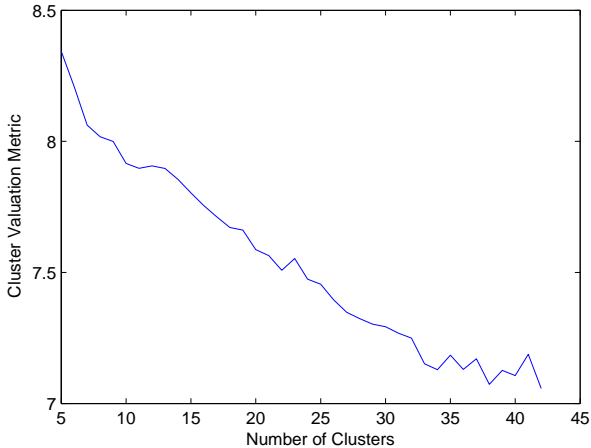}
	\caption{Cluster Valuation Metric for Given Numbers of Clusters Using $k$-means}
	\label{cluster-sizes}
\end{center}
\end{figure}

The $k$-means clustering algorithm \cite{kmeans} with Euclidean distances is
used on the transport load vectors.  BEVs and PHEVs were separated before the
clustering so that vehicles corresponding to the same base driving profile are
either all BEVs or all PHEVs.  Figure \ref{cluster-sizes} shows the cluster
valuation metric we used to determine the benefit of the number of clusters.
Each cluster has a centroid (or mean), and we measure the distance from each
point in the cluster to its centroid, called the within-cluster mean. We then
find the mean of these distances among all clusters, say the overall-cluster
mean.  Since $k$-means clustering is a randomized algorithm, we perform 50 runs
of the algorithm on different fleets of transport loads $t_j$, and take the
statistical mean of the overall-cluster mean over all the runs; the logarithm
of this value is our cluster valuation metric on the $y$-axis of Figure
\ref{cluster-sizes}, for the range of 5 to 45 clusters.  As Figure
\ref{cluster-sizes} shows, the valuation metric within-cluster means stop
decreasing at $k=37$ and no more than 43 clusters can be formed; thus, we form
37 base driving profiles.  By clustering 400 driving profiles from a similar
day (i.e., same day of the week in the same region) into 37 base driving
profiles and determining the portion of vehicles in each cluster, expected
hourly demand of a fleet of 10,000 PEVs can be estimated.

\subsection{Dual Problem of the Clustered Formulation}

We form 37 clusters of driving behaviors for a similar previous $n$-hour
period, and assign each cluster $\ell = 1, \ldots, 37$ a weight $b_\ell$. This
weight takes into account both the number of vehicles and the battery sizes of
vehicles in the group, so that the demand of each cluster is proportionally
accounted for.  Since the base driving schedule of each cluster is known for
the entire $n$-hour period, we can use them to form a static formulation as a linear
program to determine feasible charging, fueling and generating schedules,
$c_\ell, f_\ell, g_\ell \in \R^n$, while minimizing the cost and staying within
the charge cap. This formulation will act as a proxy for the aggregate 
behavior of the actual fleet:
\begin{equation}
	\begin{array}{lll}
		\underset{c_\ell, f_\ell, g_\ell, \forall \ell}{\minimize} 
			& \sum_{\ell=1}^k b_\ell \cdot
			\left(p^T c_\ell + (p^g)^Tf_\ell \right) \\
		\st & s_{\ell,1:n} = s_{\ell,1:n-1} + c_{\text{eff}} 
			\cdot c_\ell + g_{\text{eff}} \cdot g_\ell-
			(\kpm) \cdot \widehat{t_\ell} 
			& : \lambda^s_\ell \\
		 & s^g_{\ell,1:n} = s^g_{\ell,1:n-1} + f_\ell - 
		(\text{gallon/kWh}) \cdot g_\ell & : \lambda^{s^g}_\ell\\
		& t_{\ell h} = 0 \;\implies\; f_{\ell h} = 0, g_{\ell h} = 0 
			&: \eta^1_\ell, \eta^2_\ell \\
		& t_{\ell h} > 0 \; \implies \; c_{\ell h} = 0 & : 
			\eta^3_\ell \\
		& \sum_{\ell=1}^m b_\ell \cdot c_\ell \le c_{\text{cap}} & : 
			\theta \\
		& 0 \le (s_\ell,c_\ell,s^g_\ell,f_\ell,g_\ell) \le 
			(\bar{s}_\ell,\bar{c}_\ell ,\bar{s}^g_\ell,
			\bar{f}_\ell,\bar{g}_\ell) & : \nu_\ell, \gamma_\ell 
	\end{array}
	\label{lp_clst}
	\tag{CLP}
\end{equation}
where $\widehat{t_\ell} \in \R^n$ is the mean driving load of the vehicles in cluster
$\ell$, for $\ell = 1, \ldots, k$.  The variables on the far right in
\eqref{lp_clst} are the dual variables associated with each constraint, where
\[
\lambda^s_\ell, \lambda^{s^g}_\ell \in \R^{n+1}, \quad
\eta^1_\ell, \eta^2_\ell, \eta^3_\ell \in \R^n, \quad
\theta \in \R^n, 
\]
and
\[
\nu_\ell = (\nu^s_\ell, \nu^c_\ell, \nu^{s^g}_\ell, \nu^f_\ell, 
	\nu^g_\ell) \in \R^{5n+2}, \quad
\gamma_\ell = (\gamma^s_\ell, \gamma^c_\ell, \gamma^{s^g}_\ell, 
	\gamma^f_\ell,\gamma^g_\ell) \in \R^{5n+2}.
\]
For simplicity, define the vector of all variables for each cluster as 
$x_\ell = (s_\ell, c_\ell, s^g_\ell, f_\ell, g_\ell) \in \R^{5n+2}$.
The linear program \eqref{lp_clst} has the dual problem:
\begin{equation}
	\begin{array}{lll}
	\maximize & \sum_{\ell=1}^k \left( (\lambda^s_\ell)^T 
		\bmat{-s_{0\ell} \\ \widehat{t_\ell}} + 
		(\lambda^{s^g}_\ell)^T \bmat{-s^g_{0\ell} \\ 0}  + 
			\gamma_\ell^T x_\ell 
		\right) + \theta^T c_{\text{cap}}   \\
	\st & \lambda^s_{0:n-1,\ell} - \lambda^s_{1:n,\ell} + 
		\gamma^s_\ell + \nu^s_\ell = 0 & : s_\ell \\
	& p - {c_{\text{eff}}}^\ell \lambda_{1:n,\ell} - I_d^\ell 
		\cdot \eta^1_\ell - b_\ell \cdot \theta	+ \gamma^c_\ell 
		- \nu^c_\ell = 0 &: c_\ell \\
	& \lambda^{s^g}_{0:n-1}-\lambda^{s^g}_{1:n} + \gamma^{s^g} 
		- \nu^{s^g} = 0 & :s^g_\ell \\
	& p_g - \lambda^{s^g}_{1:n,\ell} - I_{nd}^\ell \eta^2_\ell 
		+ \gamma^f_\ell - \nu^f_\ell = 0 & : f_\ell\\ 	
	& - {g_{\text{eff}}}^\ell \cdot \lambda^{s}_{1:n,\ell} - 
	(\text{gallon/kWh}) \cdot \lambda^{s^g}_{1:n,\ell}-I_{nd}^\ell
		\eta^3_\ell + \gamma^g_\ell - \nu^g_\ell = 0 & : g_\ell \\
	& \nu_\ell, \gamma_\ell \geq 0
	\end{array}
	\label{dual_clst}
	\tag{DCLP}
\end{equation}
where $I_d^\ell \in \R^{n\times n}$ is the identity matrix with the columns
corresponding to the hours cluster $\ell$ is \emph{not} driving `zeroed out',
and $I_{nd}^\ell\in \R^{n\times n}$ is the identity matrix with the columns
corresponding to hours the cluster \emph{is} driving `zeroed out',
i.e., $I_d^\ell + I_{nd}^\ell = I_n$.

In contrast to the large problem \eqref{lp} with tens of thousands of vehicles,
this formulation is considerably smaller, with one schedule per cluster.  In
our testing, once the number of clusters to form was determined, forming a set
of clusters took around 0.3 seconds and solving the clustered linear program
\eqref{lp_clst} to find the corresponding dual variables can be done in about
40 seconds on a single workstation.  Moreover, the solution and dual variables
are only found once, before the dynamic portion of the algorithm assigns
charging schedules.  

\section{Constraint-Adjusted Pricing}

Our dynamic algorithm uses properties of the dual problem \eqref{dual_clst} to
create adjusted `prices' for each behavior cluster, so that allocating
charging among the hours with the lowest adjusted price value will result in a
feasible and satisfactory solution.  These new prices are called
\emph{constraint-adjusted prices}, since they use the dual variables in the
constraints of \eqref{dual_clst}.  Specifically, we use the second constraint
in \eqref{dual_clst} corresponding to the primal variable $c_\ell$, 
\[
	p - {c_{\text{eff}}}^\ell \lambda_{1:n,\ell} - I_d^\ell \cdot \eta^1_\ell 
		- b_\ell \cdot \theta	+ \gamma^c_\ell - \nu^c_\ell = 0, 
\]
where $\gamma^c, \nu^c$ are slack variables. By strict complementarity,
for each hour $h = 1, \ldots, n$, 
\begin{align*}
	\text{if }  & (p_{h} - c_{\text{eff}}^\ell \cdot \lambda_{h\ell} 
		- I_d \cdot \eta^1_{h\ell} - b_\ell \cdot \theta_h) > 0, \text{ then }
	c_{\ell h} = 0 \\
	\text{and if }  & (p_{h} - c_{\text{eff}}^\ell \cdot \lambda_{h\ell} 
		- I_d \cdot \eta^1_{h\ell} - b_\ell \cdot \theta_h) < 0, \text{ then }
	c_{\ell h} = \bar{c_{\ell}} 
\end{align*}
where $\bar{c}_{\ell}$ is the maximum charging rate of cluster $\ell$.  The
value
\[
	d_\ell := p_h - c_{\text{eff}}^\ell \cdot \lambda_{h\ell} - I_d \cdot \eta^1_{h\ell} - 
		b_\ell \cdot \theta_h
\]
can be thought of as a vector of `constraint-adjusted prices' over the next $n$
hours for cluster $\ell$ that takes into account: the charge cap, current
electricity price, current battery energy storage, and the vehicle's driving
schedule.

\begin{figure}[!htbp]
\hspace{-1in}
\begin{center}
	\includegraphics[scale=.55]{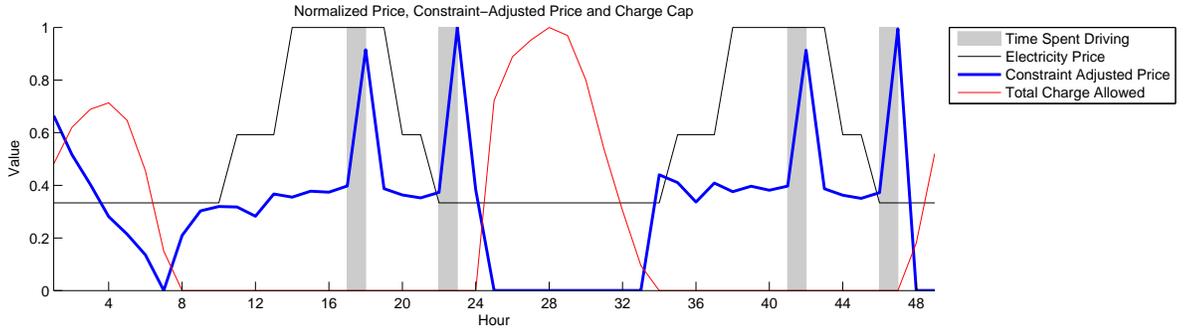}
\caption{Comparison of the Normalized Constraint-Adjusted Price vs.\ Electricity Price}
\label{CAP_compare}
\end{center}
\end{figure}

Figure \ref{CAP_compare} shows an example of the constraint-adjusted prices,
along with the electricity price, total charge allowed, and driving schedule
for a vehicle over a 2-day (or 48-hour) period, where each set of values is
normalized by its respective largest element (i.e., normalized by the
$\ell_1$-norm). Here, the `total charge allowed' is the amount of electricity
the scheduling obligation allows for the fleet to charge, i.e., it is the
electricity load without PEVs subtracted from the charge cap.  For this
example, the driver's profile is assumed to be the same on each day.  The
electricity load is based on a week in August $2011$ in the PG\&E transmission
access charge area, and the prices used are PG\&E baseline summer time-of-use
(TOU) rates.  The charge cap at each hour is assumed to be the daily peak
demand value; in broad terms, this will mean that the total charge allowed
increases as the demand decreases, and no charging is allowed during the hour
when total electricity demand reaches its peak.

Note that the constraint-adjusted price does not actually become zero; since
these are the normalized values, this plot shows that the constraint-adjusted
price becomes relatively very small when the supply is large.  Furthermore,
the constraint-adjusted price is relatively large during the first hours,
because we assume each vehicle starts with a full battery.  

The constraint-adjusted prices represent the cost to the aggregator of a
vehicle charging in a certain hour when considering the above factors, as is
shown in Figure \ref{CAP_compare}.  The price increases when the vehicle is
driving, and when the charge cap is very small.  Note that the
constraint-adjusted price still takes into account the hourly electricity
price, so the vehicle will first charge in the hours with the lowest cost, but
the additional factors are also taken into account.  A mechanism using the
constraint-adjusted prices to determine vehicle charging, as opposed to
electricity price, would assign charging schedules that satisfy both the
consumers and the market.

As each vehicle plugs in to the grid, it reports its expected driving schedule
for the next $n$ hours.  Algorithm \ref{CoAdPr} uses these constraint-adjusted
prices for the next $n$ hours to instantly determine low-cost feasible charging
schedules.  The mechanism will charge the battery (and notify the driver to
fill the tank, if necessary) the amount needed to complete a given transport
load, during the feasible hours where the constraint-adjusted price is lowest.
Since Step \ref{stp2} is executed only once and prices are constant over the
horizon, the entire charging schedule for each vehicle is determined at the
time it connects.  This is based on an assumption that electricity prices are 
constant, which may be relaxed in future work.  
\begin{algorithmic}[1]
	\begin{algorithm}
	\caption{Constraint-Adjusted Pricing}
	\label{CoAdPr}
		\STATE Cluster vehicle driving patterns from a similar day using $k$-means 
			with Euclidean distances \label{stp1}
		\STATE Solve for the dual variables of the problem \eqref{lp_clst} using
			the clusters from Step \ref{stp1} \label{stp2}
		\FOR{each vehicle that plugs in}
			\STATE Determine which cluster $\ell$ the vehicle belongs to and its
				constraint-adjusted price $d^\ell$.
			\WHILE{the vehicle needs energy}
				\IF{the vehicle is connected}
					\STATE Charge the battery the required amount during hours with the
					smallest values of $d^\ell$, while staying within the charge cap. 
				\ELSIF{the vehicle has a gasoline tank and cannot receive sufficient 
					energy from charging}
						\STATE Generate electricity from gasoline and activate
						notification to fill gasoline tank when needed
				\ENDIF
			\ENDWHILE
		\ENDFOR
	\end{algorithm}
\end{algorithmic}

Revisiting the previous example of two vehicles in three hours, the
constraint-adjusted prices for vehicle $j=1,2$ are
\[
	d_1 = \bmat{0.38 \\ 0 \\ 0.76}, \quad 
	d_2 = \bmat{0 \\ 0.76 \\ 0.07}.
\]
If the two vehicles are assigned charging schedules with constraint-adjusted
prices via Algorithm \ref{CoAdPr}, vehicle 1 will charge in hour 2 (not hour 1,
as before), allowing vehicles 2 to charge in hour 1.

\section{Comparison Algorithms}
\label{compare}

We compare Constraint-Adjusted Pricing, Algorithm \ref{CoAdPr}, to a number of
other possible charging scenarios, resulting from either standard charging
without a DR mechanism or from a DR mechanism that allocates charging
differently.  We compare three standard measurement and benefit features
\cite{ferc-DR} from each charging scenario: cost to the consumer, increase in
peak electricity demand and the total amount of energy used. 

In general, the \emph{baseline} for comparison is the estimate of electricity
usage in the absence of a DR mechanism  \cite{enernoc-baseline}. Since we are
simulating a fleet of vehicles, we can compare the results of our algorithm to
the outcome if vehicles were to charge each time they plug in for at least an hour
(Algorithm \ref{alg-sc}, Standard Charging).  We also compare our results to
two potential DR mechanisms (Algorithm \ref{lowcost}, Lowest-Cost Charging and
Algorithm \ref{pri_perc}, Relative Primal Percentages).  We briefly describe
each algorithm used for comparison below before summarizing the results.

\subsection{Standard Charging}
For a standard baseline, we simulate the outcome of vehicles automatically
charging each time they plug-in, until either the battery is full or the
connection is ended; we assume that a vehicle is plugged in to the grid each
time it is parked for at least an hour.  This is called Standard Charging, as
described in Algorithm \ref{alg-sc}, since these charging schedules are the
outcome of an absence of DR intervention.
\begin{algorithmic}[1]
	\begin{algorithm}
	\caption{Standard Charging}
	\label{alg-sc}
		\FOR{each vehicle that plugs in to the grid}
			\WHILE{the battery is not full and the vehicle is connected}
				\STATE Charge the battery at the maximum rate 
			\ENDWHILE
		\ENDFOR
	\end{algorithm}
\end{algorithmic}

\subsection{Lowest-Cost Charging}
For a comparison to DR mechanisms that are currently in place, we simulate the
algorithm where vehicles only charge when needed, in the hours with the lowest
possible electricity price, as in Algorithm \ref{lowcost}.  This is equivalent
to a user charging only in hours with the lowest electricity price (i.e., when
the price is below a given threshold) unless such charging is not sufficient,
in which case the battery is also charged in hours with a higher electricity
price (i.e., when the price is above a given threshold). As described in
\cite{enernoc-DR-types}, this is an Economic-based DR program, which is a
concept used in practice with DR resources for the Commercial and Industrial
sector.

\begin{algorithmic}[1]
	\begin{algorithm}
	\caption{Lowest-Cost Charging}
	\label{lowcost}
		\FOR{each vehicle that plugs in to the grid}
			\WHILE{the vehicle needs energy}
				\IF{the vehicle is connected}
					\STATE Charge the battery the required amount in hours with the
						lowest electricity price
				\ELSIF{the vehicle is driving and there is a gasoline tank}
					\STATE Generate electricity from gasoline and if necessary, activate
					notification to fill gasoline tank 
				\ENDIF
			\ENDWHILE
		\ENDFOR
	\end{algorithm}
\end{algorithmic}

\subsection{Relative Primal Percentages}
The last algorithm we use for comparison uses the solution of the primal
problem \eqref{lp_clst} to determine charging, generating and fueling schedules
of a fleet of PEVs.  The solution to the primal problem \eqref{lp_clst} takes
into account the charge cap, electricity price, state of the battery, and
driving schedule.  Thus, we show that generating charging schedules using the
primal solution as a guideline is not as effective as using information from
the dual problem to construct Constraint-Adjusted Pricing.

\begin{algorithmic}[1]
	\begin{algorithm}
	\caption{Relative Primal Percentages}
	\label{pri_perc}
		\STATE Cluster vehicle driving patterns from a similar day using $k$-means
			with Euclidean distances \label{stp1-2}
		\STATE Solve for the primal variables in the problem \eqref{lp_clst} using the
			clusters from Step \ref{stp1-2}
		\FOR{each vehicle that plugs in}
			\STATE Determine which cluster $\ell$ the vehicle belongs to 
			\IF{the vehicle is connected}
				\STATE Charge the battery in the same portions of the needed amount as
					cluster $\ell$ does in \eqref{lp_clst}
			\ELSIF{the vehicle is driving and there is a gasoline tank}
				\STATE Generate electricity from gasoline in the same portions of the
				needed amount as corresponding cluster $\ell$ does in \eqref{lp_clst}
				and if necessary, activate notification to fill gasoline tank
			\ENDIF
		\ENDFOR
	\end{algorithm}
\end{algorithmic}

This algorithm assigns schedules so that each vehicle depends on power from
charging, generating, and fueling in the same ratios and in the same relative
percentage in each hour as in the primal.  For example, if a cluster obtains
$\alpha_c\%$ of its energy from charging, $\alpha_g\%$ from generating, and
$\alpha_f\%$ from fueling, then so will every vehicle in this cluster.
Moreover, if $r^c_{\ell 1}, r^c_{\ell 2},\ldots r^c_{\ell n}$ are the
respective ratios of total charging cluster $\ell$ does in each hour, where
$\sum_{h=1}^n r^c_{\ell h} = 1$, then the amount of charging in each hour for
any vehicle in the cluster will be in the same ratios, within feasible bounds.
Similarly, the vehicle will generate and fuel in the ratios ${r^g_\ell, r^f_\ell
\in \R^n}$.  Algorithm \ref{pri_perc} is a DR algorithm based on the primal
solution of \eqref{lp_clst} to determine what portion of the over all charging
should be done in each hour. 

\section{Results}

The following examples show the resulting charging schedules for a fleet of
10,000 PEVs over a 5-day period.  We compare Algorithm \ref{CoAdPr} to the
algorithms in Section \ref{compare} to determine the advantages of 
Constraint-Adjusted Pricing.  We use real data, including: %
\begin{itemize}\compresslist
	\item driving patterns in urban California, obtained from the NHTS 
	dataset \cite{nhts},
	\item electricity loads in the PG\&E transmission area for the week of August
		22-28, 2011 \cite{caiso-oasis},
	\item PG\&E baseline summer time-of-use rates (electricity 
		pricing) \cite{pge-tou}, 
	\item and the mean gasoline price in the zip code 94305 on August 
		31, 2011 \cite{gas-price}.
\end{itemize}

\begin{figure}[!h]
	\begin{center}
	\includegraphics[scale=0.65]{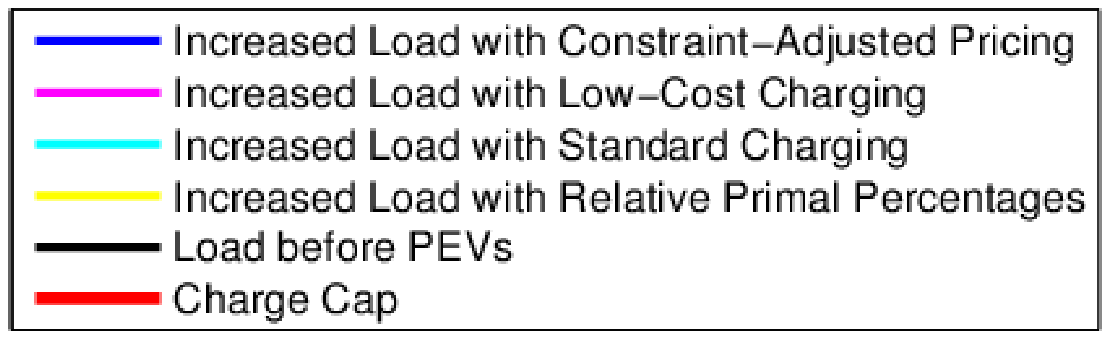}
	\end{center}
	\vspace{-.1in}

	\hspace{-1.2in}
	\includegraphics[scale=0.7]{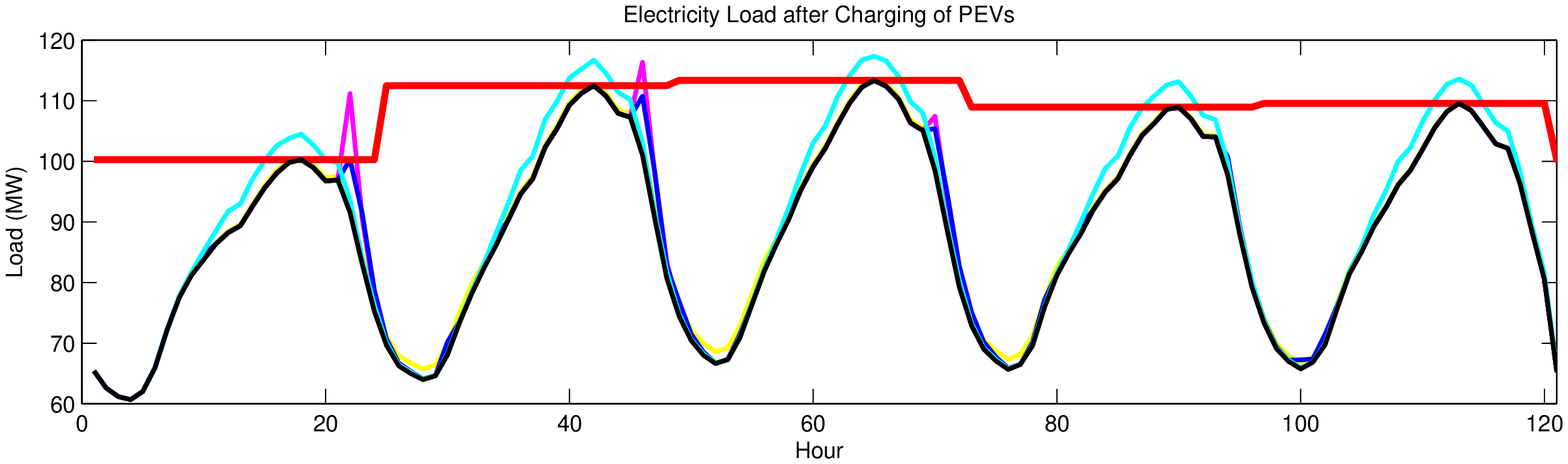}
	\vspace{-.45in}
	\begin{center}
	{\small
	\begin{tabular}{|l|c|c|c|c|}
		\hline
		& Standard &Lowest-Cost &  Relative Primal & Constraint-Adjusted  \\
		& Charging & Charging & Percentages & Pricing \\
		\hline
		Increase in Peak Power Demand (\%) & 4.2 & 10.9 & 0.02 & 0 \\
		\hline
		Total Grid Energy Fleet Used (MWh) & 255.3 & 231.5 & 231.5 & 231.5 \\
		\hline
		Total Gasoline Fleet Used (MWh) & 0.4 & 0.3 & 0.3 & 0.3 \\
		\hline
		Total Cost (\$1,000)	& 99.6 & 65.6 & 67.7 & 65.6 \\
		\hline
	\end{tabular}}
	\end{center}
	\caption{Charge cap at 100\% of the peak}
	\label{charge_cap100}
\end{figure}

\begin{figure}[!h]

	\hspace{-1.2in}
	\includegraphics[scale=0.75]{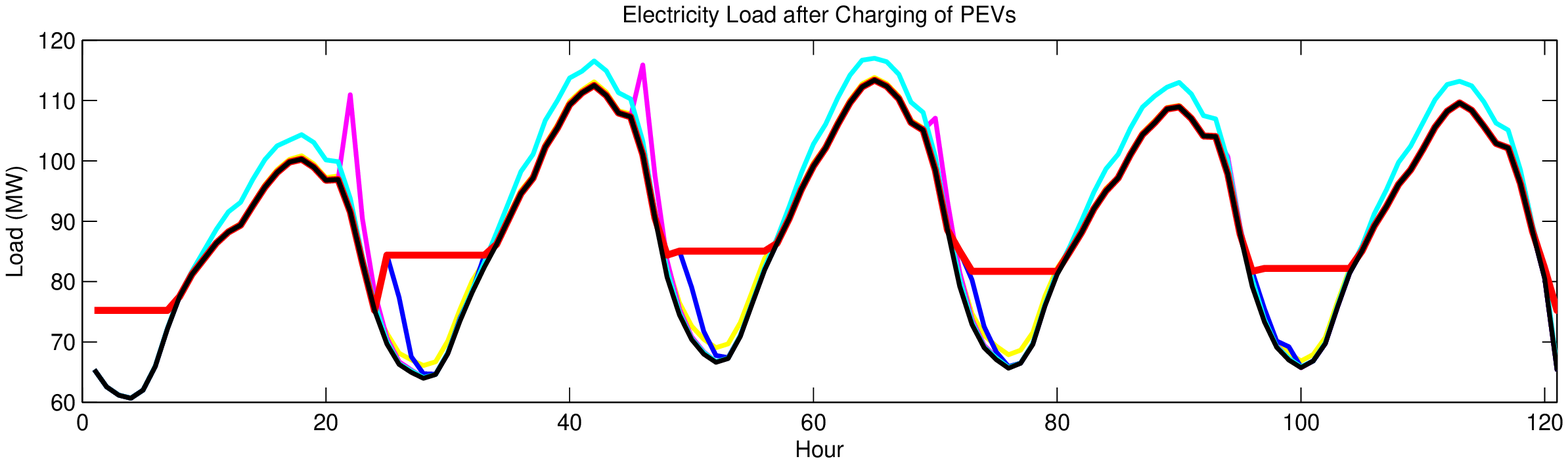}
	\vspace{-.45in}
	\begin{center}
	{\small
	\begin{tabular}{|l|c|c|c|c|}
		\hline
		& Standard & Lowest-Cost & Relative Primal & Constraint-Adjusted \\
		& Charging & Charging & Percentages & Pricing \\
		\hline
		Increase in Peak Power Demand (\%) & 4 & 10.6 & 0.5 & 0 \\
		\hline
		Total Grid Energy Fleet Used (MWh) & 254.2 & 230.7 & 230.9 & 230.7 \\
		\hline
		Total Gasoline Fleet Used (MWh) & 0.4 & 0.3 & 0.3 & 0.3 \\
		\hline
		Total Cost (\$1,000) & 97.6 & 64.0 & 66.4 & 64.0 \\
		\hline
	\end{tabular}}
	\end{center}
	\caption{Charge cap at 75\% of the peak}
	\label{charge_cap75}
\end{figure}

Figures \ref{charge_cap100} and \ref{charge_cap75} are examples that show the
charging demands of a fleet of 10,000 PEVs added to the electricity load, using
each algorithm.  Lowest-Cost Charging and Standard Charging clearly create new
peak demands, while charging the fleet with Constraint-Adjusted Pricing ensures
the load stays within the charge cap.

The statistics in Tables \ref{pwr-energy} and \ref{costs} are the mean values
taken over 250 simulations of a fleet with 10,000 vehicles for 5 days, or 120
hours.  Each simulation uses a different fleet of vehicle driving behaviors,
sampled from the NHTS data.  Note that the clusters are formed from a separate
subset of the data that is not used for testing.  We assume that there is one
PEV for every 3 households (the electricity demand is for 30,000 households)
and that each vehicle plugs in over a 12-hour period, giving its expected
driving schedule until the end of the fifth day.  We also assume each vehicle
has a full gasoline tank and battery at the beginning of the 5-day period.  The
costs listed below account for the cost of electricity and gasoline used to
drive the vehicles, i.e. the costs account for new charging done and previously
stored energy that is used for driving.

Constraint-Adjusted Pricing is the only mechanism simulated that does not
create a new peak electricity demand.  With a penetration of 33\% of households
owning a PEV, Standard Charging results in a 3.5\% increase in peak power
demand, and Lowest-Cost Charging (of Economic-based DR programs) result in a
9.8\% increase.

The total consumer cost of charging a PEV using Constraint-Adjusted Pricing is
a reduction of more than 30\% of the consumer cost from Standard Charging, and
equal to the Lowest-Cost charging.  However, Constraint-Adjusted Pricing is
significantly better than Lowest-Cost Charging in terms of peak power
increase, that is, Lowest-Cost Charging will not satisfy a market charge cap.
Moreover, our aggregated linear program \eqref{lp_clst}, which finds the
constraint-adjusted prices, was solved in under a minute on a single
workstation.  Therefore, when considering consumer costs along with an
increase in peak power demand, Constraint-Adjusted Pricing out-performs each
algorithm in this comparison.

\begin{table}[!htbp]
\begin{center}
{\small
\begin{tabular}{|l|c|c|c|c|}
	\hline
	& Standard  & Lowest-Cost & Relative Primal & Constraint-Adjusted\\
	& Charging & Charging & Percentages & Pricing  \\
	\hline
	Increase in Peak Power Demand (MW) & 3.9 & 9.8 & 1.9e-2 & 0 \\
	\hline
	Percent Increase in Peak Power Demand (\%)  & 3.5 & 9.8 & 1.7e-2 & 0 \\ 
	\hline
	Total Grid Energy Fleet Used (MWh) & 248.5 & 225.6 & 228.6 & 225.6 \\
	\hline
	Total Gasoline Fleet Used (MWh) & 0.36 & 0.34 & 0.33 & 0.34 \\
	\hline
\end{tabular}}
\caption{Power and Energy Comparisons: Mean Values}
\label{pwr-energy}
\end{center}
\end{table}

\begin{table}[!htbp]
\begin{center}
{\small
\begin{tabular}{|l|c|c|c|c|}
	\hline
	& Standard & Lowest-Cost& Relative Primal & Constraint-Adjusted  \\
	& Charging & Charging & Percentages & Pricing \\
	\hline
	Total Cost (\$1,000)& 98.5 & 64.9 & 65.7 & 64.9 \\
	\hline
	Gasoline Cost (\$1,000) & 47.7 & 44.6 & 43.7 & 44.6 \\
	\hline
	Electricity Cost (\$1,000) & 50.9 & 20.3 & 22.0 & 20.3 \\
	\hline
	Mean Cost/Mile (\$) & 0.069 & 0.030 & 0.031 & 0.030 \\
	\hline
\end{tabular}}
\caption{Consumer Cost Comparisons: Mean Values}
\label{costs}
\end{center}
\end{table}

\section{Conclusions}

We constructed a dynamic mechanism to charge a fleet of PEVs.  The allocation
of electricity can satisfy an upper limit obligation to the market and provide
each vehicle with enough energy for its transport load.  Our algorithm depends
only on the knowledge of a few hundred driving behaviors from a previous
similar day and instantly assigns charging schedules to vehicles as they
plug-in to the grid, without depending on information about the current
transportation needs and driving schedules of the other vehicles in the fleet.
We use a simple adjusted-pricing scheme to allocate charging to feasible and
satisfactory hours. Our results show that Constraint-Adjusted Pricing
out-performs each comparison algorithm when considering consumer cost, increase
in peak power demand and total energy used.   

Our future work will include incorporating vehicle-to-grid (V2G) capabilities 
of the fleet and ensuring our mechanism is robust to unexpected events.  We
also plan to determine an equilibrium price for the aggregator to use in
determining PEV charging schedules, which will not change drastically due to
slight fluctuations in demand.

\nocite{*} 
\bibliographystyle{plain} 
\bibliography{PEVs}{}

\end{document}